\documentclass[12pt]{article}
\usepackage{amssymb}
\usepackage{amsmath}
\usepackage{amscd}
\usepackage{amsthm}

\usepackage{graphicx}
\usepackage{enumerate}

\newcommand{\C}{\mathbb C}

\newcommand{\h} {\hat}
\newcommand{\mc}{\mathcal}

\newcommand{\inver}{\underset{\longleftarrow}{\lim}}

\newcommand{\fiber}{\pi^{-1}}

\newcommand{\sol}{\mathcal{S}^1}
\newcommand{\NN}{\mathcal{N}}
\newcommand{\RR}{\mathcal{R}}

\newtheorem{theorem}{Theorem}
\newtheorem{lemma}[theorem]{Lemma}
\newtheorem{corollary}[theorem]{Corollary}
\newtheorem{proposition}[theorem]{Proposition}

\newcommand{\parag}[1]{
\medskip
\noindent {\bfseries #1}
}%---paragraph
\newcommand{\Cbar}{\overline{\mathbb{C}}}
\newcommand{\Dbar}{\overline{\mathbb{D}}}
\newcommand{\R}{\ensuremath{\mathbb{R}}}
\newcommand{\D}{\ensuremath{\mathbb{D}}}
\newcommand{\Z}{\ensuremath{\mathbb{Z}}}

\newcommand{\BB}{\mathcal{B}}

\newcommand{\II}{\mathcal{I}}%

\newcommand{\XX}{\mathcal{X}}
\newcommand{\YY}{\mathcal{Y}}

\newcommand{\QQ}{\mathcal{Q}}
\newcommand{\kakko}[1]{{\left( #1 \right)}}
\newcommand{\skakko}[1]{{\left\{ #1 \right\}}}
\newcommand{\ee}{~=~}
\newcommand{\dee}{~:=~}

\newcommand{\st}{\,:\,}
\newcommand{\QED}{\hfill $\blacksquare$}
\newcommand{\al}{{\alpha}}
\newcommand{\zhat}{\hat{z}}
\newcommand{\fhat}{\hat{f}}
\newcommand{\ghat}{\hat{g}}
\newcommand{\fc}{{f_c}}
\newcommand{\invlim}{\underset{\longleftarrow}{\lim}}
\newcommand{\ql}{quadratic-like~}
\newcommand{\sminus}{-} % tentative
\newcommand{\homeo}{\approx} % tentative
\title{Topology of the regular part for infinitely renormalizable
quadratic polynomials}
\author{Carlos Cabrera and Tomoki Kawahira}

\begin{document}
\maketitle

\begin{abstract}
In this paper we describe the well studied process of
renormalization of quadratic polynomials from the point of view of
their natural extensions. In particular, we describe the topology of
the inverse limit of infinitely renormalizable quadratic polynomials
and prove that when they satisfy a-priori bounds, the topology is
rigid modulo its combinatorics.

\end{abstract}

\section{Introduction and basic theory}
The last quarter of the last century witnessed an explosion of
results concerning the quadratic family. Of particular importance
was the development of the notion of renormalization which allowed
to describe much of the dynamical richness the family posses. In
this setting, important contributions were given by the work of
several people: Feigenbaum, Douady, Hubbard, Sullivan, Yoccoz,
Lyubich and McMullen among many others.

In \cite{sullb}, Sullivan constructed a lamination by Riemann
surfaces associated to expanding maps on the circle, by using its
inverse limit. Later on in \cite{LM}, Lyubich and Minsky generalized
this construction to every rational map on the sphere. In this
setting the construction of the lamination is more involved since
the presence of critical orbits forces to consider a subset of the
inverse limit, called the regular space, provided with a finer
topology than the induced from the product topology on the inverse
limit.

Part of the program presented by Lyubich and Minsky, it was to
investigate the properties of the regular part for infinite
renormalizable polynomials.

Under the assumption of a-priori bounds, the regular part of an
infinite renormalizable polynomial $f_c$ is a lamination under the
topology induced from its inverse limit.

In this paper, we show that the topology of the regular part
determines the dynamics of $f_c$ up to combinatorial equivalence
(Main Theorem). This implies a kind of rigidity of the regular parts
associated with infinitely renormalizable maps with a priori bounds.

\parag{Outline of this paper.}
In the rest of this section we give the basic theory of the dynamics
of quadratic maps and their renormalizations. In Section 2, we
review the definition of the inverse limits and the regular parts
generated by quadratic maps. Section 3 is devoted for the statement
and the proof of the Structure Theorem (Theorem \ref{thm.blocks}),
which claims that regular parts of the persistently recurrent
infinitely renormalizable maps are decomposed into ``blocks"
according to the tree structure associated with the nest of
renormalizations. Finally in Section 4, we prove the Main Theorem
(Theorem \ref{thm.main}) stated as above.

\parag{Acknowledgements.}
We would like to thank for their
hospitality and support to the Fields institute and IMPAN (Warsaw)
were this work was carried on. We would also want to thank Misha
Lyubich for useful conversations. The first author want to thank for
the hospitality and support of IMATE (Cuernavaca). The second author
is partially supported by JSPS Grant-in-Aid for Young Scientists,
the Circle for the Promotion of Science and Engineering, and Inamori
Foundation.

\subsection{Preliminary}
We start with some basic definitions on the dynamics of quadratic
maps and inverse limits. Readers may refer \cite{DH} and \cite{L}
for dynamics of quadratic maps.

\parag{Julia and Fatou sets.}
For quadratic map $f_c(z)= z^2+c$ on the Riemann sphere $\Cbar$
with parameter $c \in \C$, the \textit{Julia set} $J(f_c)$ is
defined as the closure of the repelling periodic points of $f_c$.
Its complement $F(f_c)=\C\sminus J(f_c)$ is called the
\textit{Fatou set}. The set $K(f_c)$ of points with bounded orbit
is called the \textit{filled Julia set}. It is known that the
boundary $\partial K(f_c)$ coincides with $J(f_c)$, and that $K(f_c)$
and $J(f_c)$ are either both connected or the same Cantor set.

\parag{B\"ottcher coordinates, equipotentials, external rays.}
Throughout this paper we assume that $K(f_c)$ and $J(f_c)$ are
both connected. In this case, the set $A_c := \Cbar \sminus
K(f_c)$ is a simply connected region which consists of points
whose orbits tend to infinity. We call $A_c$ \textit{the basin of
infinity} of $f_c$. There exists a unique Riemann map $\psi_c:A_c
\rightarrow \Cbar \sminus \Dbar$, called the \textit{B\"{o}ttcher
coordinate}, that conjugates $f_c$ in $A_c$ with $w \mapsto w^2$
on $\Cbar \sminus \Dbar$ and $\psi_c(z)/z \to 1 ~(z \to \infty)$.
For $r
>1$, $E_c(r):= \psi_c^{-1}( \skakko{w \in \C \st |w|=r})$ is
called the \textit{equipotential curve} of level $r$. For $\theta
\in \R/\Z$, $R_c(\theta):= \psi_c^{-1}(\skakko{w \in \C \sminus \Dbar \st
\arg w =\theta})$ is called the \textit{external ray} of angle
$\theta$.

\parag{Ray portrait.}
Let $O=\skakko{p_1, \ldots, p_m}$ be a repelling cycle of $f_c$.
There are finitely many angles of external rays landing at each
$p_i$, which we denote by $\Theta(p_i)$. It is a fact due to Douady
and Hubbard \cite{DH} that $\Theta(p_i)$ is a set of rational
numbers. The collection $\mathrm{rp}(O)=\{\Theta(p_1), \ldots,
\Theta(p_m) \}$ is called the \textit{ray portrait} of $O$. A ray
portrait is called \textit{non trivial}, if there are at least two
rays landing at every point in $O$. A non trivial ray portrait
determines a region in the parameter space with a leading hyperbolic
component. In this way, every non trivial ray portrait determines a
unique superattracting parameter by taking the center of its leading
hyperbolic component. (See Milnor's \cite{Mext})

\parag{Superattracting quadratic maps.}
Quadratic maps $f_s(z)=z^2+s$ whose critical orbit is periodic form
an important class of quadratic maps, called
\textit{superattracting} quadratic maps. Let $\skakko{\al_s(1),
\ldots, \al_s(m)=0}$ denote the critical cycle with
$f_s(\al_s(i))=\al_s(i+1)$, where we take the indexes modulo $m$.
Then it is known that the connected component $D_s$ of the Fatou set
(``\textit{Fatou component}") with $\al_s(m)=0 \in D_s$ is a Jordan
domain with dynamics $f_s^m: \Bar{D_s} \to \Bar{D_s}$ conjugate to
$f_0: \Dbar \to \Dbar$. Let $\Psi_s: \Bar{D_s} \to \Dbar$ be this
conjugacy, which we also call a B\"{o}ttcher coordinate for
$\Bar{D_s}$. The \textit{internal equipotential} $I_s(r)$ of level
$r < 1$ is defined by $\Psi_s^{-1}(\skakko{|w|=r})$. We also denote
$\Psi_s^{-1}(\skakko{|w|< r})$ by $D_s(r)$.

The pull-back of $1 \in \Dbar$ by $\Psi_s$ in $\partial D_s$ is a repelling
periodic point with period $\le m$. Let $O_s$ be its cycle which
is on the boundary of $\bigcup_{1 \le i \le m} f_s^i(D_s)$. We say
the ray portrait $\mathrm{rp}(O_s)$ is the \textit{characteristic
ray portrait} of $f_s$. In fact, superattracting $f_s$ is uniquely
identified by $\mathrm{rp}(O_s)$. (See Milnor's \cite{Mext}).

\parag{Quadratic-like maps.}
Let $U$ and $V$ be topological disks in $\C$ with $U$ compactly
contained in $V$. A \textit{quadratic-like map} $g:U \to V$ is a
proper holomorphic map of degree two. The \textit{filled Julia set}
is defined by $K(g):=\bigcap_{n \ge 1} g^{-n}(V)$. Throughout this
paper we assume that any \ql map $g:U \to V$ has a connected $K(g)$.
Its \textit{Julia set} $J(g)$ is the boundary of $K(g)$. The
\textit{postcritical set} $P(g)$ is the closure of the forward orbit
of the critical point of $g$, since $K(g)$ is connected we have
$P(g)\subset U$.

By the Douady-Hubbard straightening theorem \cite{DH}, there exists
a unique $c=c(g) \in \C$ and a quasiconformal map $h:V \to V'$ such
that $h$ conjugates $g:U \to V$ to $f_c: U' \to V'$ where
$U'=h(U)=\fc^{-1}(V')$ and $\bar{\partial} h =0$ a.e. on $K(g)$. The
quadratic map $f_c$ is called the \textit{straightening} of $g$ and
$h$ is called a \textit{straightening map}. Though such an $h$ is
not uniquely determined, we always assume that any \ql map $g$ is
accompanied by one fixed straightening map $h=h_g$.

One can easily check that there exists an $r_g>1$ such that if $1<r
\le r_g$ and $\theta \in \R/\Z$, the pulled-back equipotentials and
external rays
$$
E_g(r) \dee h^{-1}(E_c(r)) ~~\text{and}~~ R_g(\theta) \dee
h^{-1}(\skakko{\rho e^{2 \pi i \theta} \st 1<\rho \le r_g})
$$
are defined. For the straightening $f_c$ of $g$, there exists a
repelling or parabolic fixed point $\beta(f_c) \in K(f_c)$ which is
the landing point of the external ray $R_c(0)$. Note that
$\beta(f_c)$ is repelling unless $c=1/4$.  We set
$\beta(g):=h^{-1}(\beta(\fc))$ and call it the $\beta$-\textit{fixed
point} of $g$.

\parag{Renormalization of quadratic maps}
A quadratic-like map $g:U \rightarrow V$ is said to be
\textit{renormalizable}, if there exist a number $m > 1$, called the
\textit{order of renormalization}, and two open sets $U_1 \subset U$
and $V_1 \subset V$ containing the critical point of $g$, such that
$g_1=g^m|~ U_1\rightarrow V_1$, called a pre-renormalization of $g$,
is again a quadratic-like map with connected Julia set $K(g_1)$.
We say $g_1: U_1\rightarrow V_1$ is a \textit{renormalization} of $g:U\to V$.
We call $K_1:=K(g_1)$, $g(K_1)$, $\cdots$, $g^{m-1}(K_1)$ the \textit{little Julia sets}.
We also assume that $m$ is the minimal order with this
property and that $K_1$ has the following property:
\textit{For any $1 \le i<j \le m$, $g^i(K_1) \cap g^j(K_1)$ is empty
or just one point that separates neither $g^i(K_1)$ nor $g^j(K_1)$}.
Such a renormalization is called \textit{simple} or
\textit{non-crossing}. See \cite{Mc1} or \cite{Mil.loc.con} for
examples of \text{crossing} renormalizations.

%For $1 \le i \le m$, set $K_1(i):=g^i(K_1)$. Its boundary
%$J_1(i):=\partial K_1(i)$ is called a \textit{little Julia set}. On
%what follows, we will keep in mind that every structure involved in
%$K_1=K_1(m)$ is  transferred to $K_1(i)$ by $g^i$.

\parag{Infinitely renormalizable maps.}
In this paper we only deal with \ql maps which are restrictions of
some iterated quadratic map. For any quadratic map $f_c$ and any
$r>1$, $f_c|~ U_c(r) \mapsto U_c(r^2)$ is a quadratic map. Set
$g_0=f_c$, $U_0:=U_c(r)$ and $V_0:=U_c(r^2)$. We say $f_c$ is
\textit{infinitely renormalizable} if there is an infinite sequence
of numbers $p_0=1 <p_1< p_2< \cdots $ and two sequences of open sets
$\{U_n\}$ and $\{V_n\}$ such that each $g_{n}=f_c^{p_n}:U_{n} \to
V_{n}$ is a \ql map, with the property that $g_{n+1}$ is a
pre-renormalization of $g_{n}$ of order $m_{n}:=p_{n+1}/p_{n}>1$.
See \cite{LyuI} for more details. The index $n$ of $g_n$ is called
the \textit{level} of renormalization.

\parag{Combinatorics of renormalizable maps.}
For a complete exposition of combinatorics of renormalizable maps we
refer to the work of Lyubich in \cite{Lyucomb} and \cite{LyuI}. From
now on $f_c$ will denote an infinitely renormalizable quadratic map
and $\skakko{g_n:U_n \to V_n}$ be its associated sequence of \ql
maps as above. In order to describe the combinatorics of $f_c$,
first we observe that the orbit of the $\beta$-fixed point of
$g_{n+1}$ forms a repelling cycle $O_{n}$ of $g_n$. Since every
$g_n$ has a unique straightening $f_{c_n}$ with $c_n=c(g_n)$ by the
straightening map $h_n$, $h_n(O_n)$ is also a repelling cycle of
$f_{c_n}$ with at least 2 external rays landing at each point in
$h_n(O_n)$, hence its ray portrait $\mathrm{rp}(h_n(O_n))$ is
non-trivial. Since every non-trivial ray portrait determines a
unique superattracting quadratic map, the $n$-level of
renormalization induces a unique superattracting map
$f_{s_n}(z)=z^2+s_n$ with characteristic ray portrait
$\mathrm{rp}(h_n(O_n))$. We call the infinite sequence of
superattracting parameters $\{s_0, s_1, s_2, \ldots \}$ the
\textit{combinatorics} of $f_c$, note that the period of the
critical point of $f_{s_n}$ is equal to $m_n$.

We say that $f_c$ has \textit{bounded combinatorics} if the sequence
$\{m_n\}$ is bounded.
The polynomial $f_c$ is said to have
\textit{a-priori bounds} if there exist $\epsilon>0$, independent of
$n$, such that $mod(V_n\setminus U_n)>\epsilon$. The map $f_c$ is
called \textit{Feigenbaum} if it has \textit{a priori} bounds and
bounded combinatorics.

\section{Inverse limits and regular parts}
In the theory of dynamical systems we use the technique of the inverse limit to construct an invertible dynamics out of non-invertible dynamics. 
In this section we give some inverse limits associated with quadratic dynamics used in this paper. 
We also define the \textit{regular parts}, which is analytically well-behaved parts of the inverse limits, according to \cite{LM}. 
Readers may refer \cite{LM} and \cite{KL} where more details on the objects defined here are given.

\subsection{Inverse limits and solenoidal cones}

\parag{Inverse Limits.}
Consider $\{f_{-n}:X_{-n} \rightarrow X_{-n+1}\}_{n=1}^\infty$,
 a sequence of $d$-to-$1$ branched covering maps on the
 manifolds $X_{-n}$ with the same dimension.
 The \textit{inverse limit} of this sequence is defined as
$$
\inver (f_{-n},X_{-n}) \dee
\{
\hat{x}=(x_{0},x_{-1}, x_{-2} \ldots )
\in \prod_{n \ge 0} X_{-n} \st f_{-n}(x_{-n})=x_{-n+1}
\}.
$$
The space $\inver(f_{-n},X_{-n})$ has a \textit{natural topology}
which is induced from the product topology in $\prod X_{-n}$. The
projection $\pi: \inver(f_{-n},X_{-n}) \to X_0$ is defined by
$\pi(\hat{x}):=x_0$.

\parag{Example 1: Natural extensions of quadratic maps.}
When all the pairs $(f_{-n}, X_{-n})$ coincide with the quadratic
$(f_c, \Cbar)$, following Lyubich and Minsky \cite{LM}, we will
denote $\inver(f_c,\bar{\C})$ by $\NN_c$. The set $\NN_c$ is called
the \textit{natural extension} of $f_c$.
 In this case we denote the projection by $\pi_c: \NN_c \to \Cbar$.
 %We will mainly deal with proper subsets of such an $\NN_c$.
There is a natural homeomorphic action $\fhat_c: \NN_c \to \NN_c$
given by $\fhat_c(z_0, z_{-1}, \ldots) := (f_c(z_0),z_0, z_{-1},
\ldots)$.

Let $X$ be a forward invariant set, by $\h{X}$ we will denote the
\textit{invariant lift} of $X$, that is the set of $\h{z}\in
\mc{N}_c$ such that all coordinates of $\h{z}$ belong to $X$. 
In particular, $\hat{\infty}=(\infty, \infty, \ldots)$. 

The natural extension is not so artificial than it appears. 
For example, it is known
that if $f_c$ is hyperbolic, then $\fhat_c$ acting on $\NN_c \sminus
\skakko{\hat{\infty}}$ is topologically conjugate to a H\'enon map
of the form $(z,w) \mapsto (z^2+c-aw, z)$ with $|a| \ll 1$ acting on
the backward Julia set $J^-$. See \cite{HOII} for more details.

\if0
A point $\h{z} \in \mc{N}_c$ is called \textit{regular} if there is a
neighborhood $U_0$ of $z_0$ such that the pull back of $U_0$ along
$\h{z}$ is eventually univalent. The \textit{regular space}
$\mc{R}_f$ is the set of regular points in $\mc{N}_f$. Let
$\mc{I}_f$ denote the set of irregular points of $f$.
\fi

\parag{Example 2: Dyadic solenoid and solenoidal cones.}
A well-known example of an inverse limit is the \textit{dyadic
solenoid} $\sol := \invlim(f_0, \mathbb{S}^1)$, where $f_0(z)=z^2$
and $\mathbb{S}^1$ is the unit circle in $\C$. The dyadic solenoid
is a connected set but is not path-connected. Any space homeomorphic
to $\inver(f_0,\Cbar \sminus \Dbar)$ will be called a
\textit{solenoidal cone}. For $f_c$ with connected $K(f_c)$, we have
an important example of a solenoidal cone $\hat{A}_c :=\inver(f_c,
A_c)$ in $\NN_c$ by looking at $\inver(f_0,\Cbar \sminus \Dbar)$
through the inverse B\" ottcher coordinate $\psi_c^{-1}$. 
More precisely, the set $\hat{A}_c$ is given by $\hat{\psi}_c^{-1}(\inver(f_0,\Cbar \sminus \Dbar))$ where $\hat{\psi}_c^{-1}: (z_0, z_{-1}, \ldots) \mapsto (\psi_c^{-1}(z_0), \psi_c^{-1}(z_{-1}), \ldots).$  
Then $\hat{A}_c \sminus \{\hat{\infty}\}$ is foliated by sets of the form
$\mc{S}_r:=\pi_c^{-1}(E_c(r))$ with $r>1$.
Each of $\mc\pi_c^{-1}(E_c(r))$ is homeomorphic to the dyadic
solenoid; in fact, the map $\phi_r:\mc{S}_r\rightarrow \sol$ given
by $\phi_r:(z_0,z_1,...) \mapsto (z_0/r,z_1/r^{1/2},...)$ is a canonical
homeomorphism. We call such $\mc{S}_r$ a \textit{solenoidal
equipotential}.

\if0
To give a few more examples of solenoidal cones, let us introduce
some notation: Let $g:U\rightarrow V$ be a proper holomorphic map,
we might allow here $U=V$, by $\invlim(g, V)$ we denote the inverse
limit for the sequence $\skakko{g: g^{-n}(V) \to g^{-n+1}(V)}_{n \ge
1}$. Let us remark that even in the cases where $g$ is defined
outside $U$, when taking preimages we will take all branches of the
inverse of $g$ satisfying $g^{-n}(V)\subset U$.
\fi

%If $g^{-1}(X) \subset X$, there is a natural action $\ghat^{-1}$ on $\invlim(g, X)$ given by $\ghat^{-1}(x_0, x_{-1}, \ldots) := (x_{-1}, x_{-2}, \ldots)$. Note that $\ghat$ is not defined on $\invlim(g, X)$ unless $X$ is a forward invariant set of $g$.

Let us give a few more examples of solenoidal cones. 
For any $r >1$, set $\D_r := \skakko{|z|<r}$.
We denote the inverse limits associated with the backward dynamics 
$$
\cdots~\rightarrow~
\Cbar \sminus f_0^{-2}(\Dbar_{r}) 
~\rightarrow~
\Cbar \sminus f_0^{-1}(\Dbar_{r}) 
~\rightarrow~ \Cbar \sminus \Dbar_r
$$
of $f_0(z)=z^2$ by $\hat{A}_0(r)$. 
This is a sub-solenoidal cone compactly contained in $\hat{A}_0 \subset \NN_0$. 
Similarly, we have a sub-solenoidal cone of $\hat{A}_c \subset \NN_c$ given by $\hat{A}_c(r):=\hat{\psi}_c^{-1}(\hat{A}_0(r))$. 
Note that the
boundary of $\hat{A}_c(r)$ in $\NN_c$ is $\mc{S}_r$. We call the
union $\hat{A}_c(r) \cup \mc{S}_r$ a \textit{compact solenoidal cone
at infinity}.

Let $f_s$ be a superattracting quadratic map as in the preceding section. 
For all $r<1$, the inverse limit given by the backward dynamics 
$$
\cdots~\rightarrow~
D_s(r^{1/4}) 
~\rightarrow~
D_s(r^{1/2})
~\rightarrow~ D_s(r) 
$$
of $f_s^{m}$ is also a solenoidal cone. We denote it by $\invlim(f_s^m, D_s(r))$. We may consider $\invlim(f_s^m, D_s(r))$ as a subset of $\NN_s$ by the following embedding map: For $(x_0, x_{-1}, \ldots) \in
\invlim(f_s^m, D_s(r))$, we define $\iota:(x_0, x_{-1}, \ldots)
\mapsto (y_0, y_{-1}, \ldots) \in \NN_s$ so that $x_{-k} = y_{-mk}$
for all $k \ge 0$. Then $\hat{D_s}(r):=\iota(\invlim(f_s^m,
D_s(r)))$ is a solenoidal cone in $\NN_s$. Note that $\partial
\hat{D_s}(r)$ is a proper subset of $\pi_s^{-1}(I_s(r))$ unless
$s=0$. 
Now $\hat{D}_s(r)$, $\fhat_s(\hat{D}_s(r))$, $\ldots$, $\fhat_s^{m-1}(\hat{D}_s(r))$ are disjoint solenoidal cones in $\NN_s$. We
also call the closures of these $m$ solenoidal cones in $\mc{N}_s$ \textit{compact solenoidal cones at the critical orbit}.

\parag{Quadratic-like inverse limits.}
Let $g:U\rightarrow V$ be a proper holomorphic map,
we might allow here $U=V$, by $\invlim(g, V)$ we denote the inverse
limit for the sequence 
$$
\cdots~\rightarrow~
g^{-2}(V)~\rightarrow~
g^{-1}(V)
~\rightarrow~ V
$$
Let us remark that even in the cases where $g$ is defined
outside $U$, when taking preimages we will take all branches of the
inverse of $g$ satisfying $g^{-n}(V)\subset U$.

Here we show the following fact on the relation between inverse limits of \ql maps and its straightening: 

\begin{proposition}\label{prop_nat_ext}
Let $g:U \to V$ be a \ql map with straightening $f_c(z)=z^2+c$. Then
the inverse limit $\invlim(g,V)$ is homeomorphic to $\NN_c$ with a
compact solenoidal cone at infinity removed.
\end{proposition}

\parag{Proof.} Set
$$
U_g(r) \dee K(g) \cup \bigcup_{1<\rho <r} E_g(\rho)
~~\text{and}~~
U_c(r) \dee K(\fc) \cup \bigcup_{1<\rho <r} E_c(\rho)
$$
for $1<r<r_g$. Then $U_g(r) \Subset V$ and $g:U_g(\sqrt{r}) \to
U_g(r)$ is a \ql map which is quasiconformally conjugate to $\fc
:U_c(\sqrt{r}) \to U_c(r)$ by straightening map $h$. Thus
$\invlim(g, U_g(r))$ is homeomorphic to $\invlim(\fc, U_c(r))$,
which is $\NN_c$ with a compact solenoidal cone $\overline{\hat{A}_c(r)}$ removed.

Now it is enough to check that the original $\invlim(g, V)$ is
homeomorphic to its subset $\invlim(g, U_g(r))$. But this follows
from the fact that $g: U \sminus U_g(\sqrt{r}) \to V \sminus U_g(r)$
is a double covering between annuli and $\invlim(g, V) \sminus
\invlim(g, U_g(r))$ is homotopic to the boundary of $\invlim(g,
U_g(r))$.\QED

\parag{Remark.}
In fact, the homeomorphism is given by a leafwise quasiconformal map
on their regular parts.

\subsection{Regular parts and infinitely renormalizable maps}

\parag{Regular parts of quadratic natural extensions.}
Let $f_c$ be a quadratic map. A point $\zhat=(z_0, z_{-1}, \ldots)$
in the natural extension $\NN_c=\invlim(f_c, \Cbar)$ is
\textit{regular} if there is a neighborhood $U_0$ of $z_0$ such that
the pull-back of $U_0$ along $\zhat$ is eventually univalent. The
\textit{regular part}(or \textit{regular leaf space})
$\RR_{f_c}=\RR_c$ is the set of regular points in $\NN_c$. Let
$\II_{f_c}=\II_c$ denote the set of irregular points.

The regular parts are analytically well-behaved parts of the natural
extensions. For example, it is known that all path-connected
components (``\textit{leaves}") of $\RR_c$ are isomorphic to $\C$ or
$\D$. Moreover, $\fhat_c$ sends leaves to leaves isomorphically.
However, most of such leaves are wildly foliated in the natural
extension, indeed dense in $\NN_c$. See \cite[\S 3]{LM} for more
details.

\parag{Example: Regular part of superattracting maps.}
A fundamental example of regular parts are given by superattracting
quadratic maps. Let $f_s$ be a superattracting quadratic map with
superattracting cycle $\skakko{\al_s(1), \ldots, \al_s(m)=0}$ as in
the previous section. Under the homeomorphic action $\fhat_s:\NN_s
\to \NN_s$, the points
$\hat{\al}_s(i):=(\al_s(i),\al_s(i-1),\al_s(i-2), \ldots )$ form a
cycle of period $m$. In this case, the set $\II_s$ of irregular
points consists of $\skakko{\hat{\infty}, \hat{\al}_s(1), \ldots
,\hat{\al}_s(m)}$. Thus the regular part $\RR_s$ is $\NN_s$ minus
these $m+1$ irregular points. Moreover, it is known that $\RR_s$ is
a Riemann surface lamination with all leaves isomorphic to $\C$.

\parag{Regular part of infinitely renormalizable maps.}
We will need the following fact, due to Kaimainovich and Lyubich,
about the topology of inverse limits of quadratic polynomials with
a-priori bounds. The proof can be found in \cite{KL}.

\begin{theorem}[Kaimainovich-Lyubich]\label{KL.loc.comp} If $f_c$
has a-priori bounds, then $\mc{R}_{c}$ is a locally compact Riemann
surface lamination, whose leaves are conformally isomorphic to
planes.\end{theorem}

\parag{Persistent recurrence.}
A quadratic polynomial $f_c:\C \rightarrow \C$ (regarded as a special case of the \ql maps) is called
\textit{persistently recurrent} if $\widehat{P(f_c)}\subset\mc{I}_c$.
Equivalently, for any neighborhood $U_0$ of $z_0 \in P(f_c)$ and any
backward orbit $\zhat = (z_0, z_{-1}, \ldots)$, pull-backs of $U_0$
along $z_0$ contains the critical point $z=0$. 
Let $f_c$ be a quadratic polynomial with \textit{a priori} bounds. If $K_n$
denotes the little Julia set of the $n$ pre-renormalization, it
follows that the postcritical set
$$P(f_c)=\bigcap_{n\geq 0}\bigcup_{j\geq 0} f_c^j(K_n)$$ is
homeomorphic to a Cantor set. Moreover, the map $f_c$ restricted to
$P(f_c)$ acts as a minimal $\mathbb{Z}$-action. 
See McMullen's \cite[Theorems 9.4]{Mc1} and the example below. 
It follows that every $f_c$ with a-priori bounds is persistently recurrent.

Hence the set of irregular points in $\invlim(f_c, \C)$ is
$\widehat{P(f_c)}$ and the projection $\pi$ restricted to
$\widehat{P(f_c)}$ is a homeomorphism over $P(f_c)$. So, we have the
following:

\begin{lemma}\label{feig.irreg} If $f_c$ is a quadratic polynomial
with a-priori bounds, then the irregular part $\mc{I}_c$ is
homeomorphic to a Cantor set together with the isolated point
$\h{\infty}$.
\end{lemma}

Let us mention that the concept of a-priori bounds is related to the
following notion of robustness due to McMullen.

\parag{Example.}
An infinitely renormalizable quadratic map $f_c$ is called
\textit{robust} if for any arbitrarily large $N>0$, there exist a
level $n>N$ of renormalization and an annulus in $\C \sminus
P(f_c)$ with definite modulus such that the annulus separates
$J(g_n)$ and $P(f_c) \sminus J(g_n)$. 
(Thus it mildly generalizes a priori bounds.)
If $f_c$ is robust, the $\omega$-limit set $\omega(c)$ of $c$ coincides with $P(f_c)$ which is a Cantor set and
the action of $f_c$ on $\omega(c)$ is homeomorphic and minimal. Thus
robust $f_c$ is also persistently recurrent. The most important
property induced by robustness is that $J(f_c)$ carries no invariant
line field, thus $f_c$ is quasiconformally rigid \cite[Theorems 1.7]{Mc1}.

\section{Structure Theorem}
In this section we will show that the natural extensions of infinitely renormalizable quadratic maps can be decomposed into
``blocks" which are given by combinatorics determined by the
sequence of renormalization.

\parag{Blocks for superattracting maps.}
We first define the blocks associated with supperattracting
quadratic maps. Let $s$ be a superattracting parameter as in Section
1, with a super attracting cycle of period $m$. For a fixed $r>1$,
we set
$$
\BB_s \dee \NN_s \sminus \overline{\hat{A}_s(r)}
\cup \bigsqcup_{i=0}^{m-1} \fhat_s^{i}(\overline{\hat{D}_s(1/r)})
$$
and call it a \textit{block} associated with $f_s$. That is, $\BB_s$
is the natural extension with compact solenoidal cones at each of
the irregular points removed. Note that $\BB_s$ is an open set and
has $m+1$ boundary components which are all solenoidal
equipotentials.

By the main result of \cite{cab} or Theorem \ref{comb.class}, if
there exists an orientation preserving homeomorphism between $\BB_s$
and $\BB_{s'}$ for some superattracting parameters $s$ and $s'$,
then $s=s'$. Thus the blocks associated with superatracting maps are
``rigid" in this sense.

In addition, we also define
$$
\QQ_s \dee \NN_s \sminus \skakko{\h{\infty}} \cup
\bigsqcup_{i=0}^{m-1} \fhat_s^{i}(\overline{\hat{D}_s(1/r)})
$$
for later use.

\parag{Structure Theorem for infinitely renormalizable maps.}
For infinitely renormalizable $f_c$ which is persistently
recurrent, it is known that $\RR_c$ is a Riemann surface
lamination with leaves isomorphic to $\C$ (Theorem \ref{KL.loc.comp}
\cite[Corollary 3.21]{KL}).
In addition, we will establish:

\begin{theorem}[{\bf Structure Theorem}]\label{thm.blocks}
Let $\fc$ be a persistently recurrent infinitely renormalizable map
and $\skakko{g_n=f_c^{p_n}|~U_n \to V_n}_{n \ge 0}$ be the
associated sequence of renormalizations with combinatorics
$\skakko{s_0, s_1, \ldots}$. Then there exist disjoint open subsets
$\BB_0, \BB_1, \ldots $ of $\NN_c$ such that:
\begin{enumerate}
\item For $n=0$, the set $\BB_0$ is homeomorphic to $\QQ_{s_0}$.
Moreover, the union $\BB_0 \cup \skakko{\h{\infty}}$
forms a neighborhood of $\h{\infty}$ with $m_0=p_1/p_0$
boundary components which are all homeomorphic to the dyadic solenoid.
\item For each $n \ge 1$, the set $\BB_n$ is
homeomorphic to $\BB_{s_n}$.
Moreover, $\BB_n$ has
$m_n+1$ (where $m_n=p_{n+1}/p_n$)
boundary components which are all homeomorphic to the dyadic solenoid.
\item For any $n \ge 1$ and $1 \le i <j \le p_n$, the sets
$\fhat_c^{i}(\BB_{n})$ and $\fhat_c^{i}(\BB_{n})$ are disjoint.
\item For $0 \le n < n'$, the closures $\overline{\BB_n}$
and $\overline{\BB_{n'}}$ intersects iff $n'=n+1$. In this case, for
all $1 \le i \le m_n$ the closures $\overline{\fhat_c^{p_n
i}(\BB_{n+1})}$ and $\overline{\BB_n}$ share just one of their
solenoidal boundary components.

\item The set
$$
\bigsqcup_{n=0}^{\infty} \bigsqcup_{i=0}^{p_n-1} \fhat_c^{i}(\BB_{n})
$$
is equal to the regular part $\RR_c$.
\item The original natural extension is given by
$$
\NN_c \ee \RR_c \sqcup \widehat{P(f_c)} \sqcup \skakko{\h{\infty}}.
$$
\end{enumerate}
\end{theorem}

By 3. and 4. above, the regular part $\RR_c$ of $f_c$ has a (locally
finite) tree structure given by configuration of blocks
$\{\fhat_c^i(\BB_n)\}$. More precisely, we join vertices
``$\fhat_c^i(\BB_n)$" and ``$\fhat_c^{i'}(\BB_{n'})$" by a segment
if they share one of their boundary component homeomorphic to
the dyadic solenoid. Then we have a \textit{configuration tree}
associated to $f_c$. Notice that, by construction, the $n$-th level
of the configuration tree of $f_c$ is a subset of the regular part
of $\mc{R}_c$. However, we do not know in general (i.e., without persistent recurrence) whether every
regular point belongs to some level of the configuration tree
associated to $f_c$.

Let us remark that the statement of Theorem \ref{thm.blocks} is
quite topological. For instance, the block $\BB_n$ which we will
construct may not be an invariant set of $f_c^{p_n}$. In the next
section, however, we will see that the topology of $\RR_c$ given by
such blocks determines the original dynamics modulo combinatorial
equivalence.

\parag{Note.}
The original motivation of this paper was to give answers to some
problems by Lyubich and Minsky \cite[\S 10]{LM}. In Problem 6, in
particular, they asked whether the hyperbolic 3-lamination
$\mc{H}_{c_0}$ and its quotient lamination $\mc{M}_{c_0}$ (they are
analogues of the hyperbolic 3-space and the quotient orbifold of a
Kleinian group) associated with the Feigenbaum parameter $c_0$,
which is the parameter of the infinitely renormalizable map
$f_{c_0}$ with combinatorics $\skakko{-1,-1, \ldots}$, reflects the
sequential bifurcation process from $f_0(z)=z^2$. In this case
$f_{c_0}$ is persistently recurrent and $\mc{H}_{c_0}$ is
constructed out of the regular part $\mc{R}_{c_0}$. Thus the
topology of $\mc{H}_{c_0}$ strongly reflects the tree structure
described by the Structure Theorem. However, since the block
decomposition of $\mc{R}_{c_0}$ is not invariant under the dynamics,
we cannot say much about the topology of the quotient lamination
$\mc{M}_{c_0}$.

A possible direction is to get $\mc{M}_{c_0}$ by a limiting process
of finitely many parabolic bifurcations. In fact, if superattracting
(or parabolic) $f_s$ is given by finitely many parabolic
bifurcations (and degenerations) from $f_0$, then the topologies of
$\mc{R}_{s}$, $\mc{H}_{s}$ and $\mc{M}_{s}$ are described in detail
\cite{Tom, Tom2}.

\subsection{Proof of the Structure Theorem}

To simplify the proof of the Structure Theorem, we first state the
main step of the proof in a proposition.

Let us start with a slightly general setting of renormalizable \ql
maps. Let $g:U \to V$ be an infinitely renormalizable \ql map with a
(simple) renormalization $g_1=g^m|~U_1 \to V_1$. Here we have to
keep in mind that we actually consider the case of $g=g_n$ and
$g_1=g_{n+1}$. But the argument also works when $g=g_0$ and $g=g_n$.
In general we do not have $V_1 \subset U$. However, we may modify $U
\Subset V$ and $U_1 \Subset V_1$ as follows: For arbitrarily fixed
$1<r<r_g$, we replace $U$ and $V$ by $U:=U_g(\sqrt{r})$ and
$V:=U_g(r)$ as in the proof of Proposition \ref{prop_nat_ext}. Note
that if we choose $r$ sufficiently close to 1 then the boundary of
$V$ is arbitrarily close to $K(g)$. Next we replace $U_1$ and $V_1$
by $U_1:=U_{g_1}(\sqrt{r_1})$ and $V:=U_{g_1}(r_1)$ with $r_1$
slightly larger than 1 so that
$$
U_1 ~\Subset~ V_1 ~\Subset~ g^{-m}(V) ~\Subset~ U ~\Subset~ V.
$$
Here the condition $V_1 \Subset g^{-m}(V)$ guarantees that the map
$g^i|_{V_1}$ makes sense and $g^i(V_1) \subset V$ for all $1 \le i
\le m$.

There exists a unique superattracting $f_s$ whose characteristic ray
portrait $\mathrm{rp}(O_s)$ is given by the cyclic orbit of
$\beta(g_1)$ by $g$. The proposition will state the relation between
$\XX:=\invlim(g,V)$, $\XX_1:=\invlim(g_1, V_1)$, and the block
$\BB_s$ associated with $f_s$ in a modified form.

Let us consider a natural embedding $\iota: \XX_1 \to \XX$ as
follows: For $\hat{x}=(x_0, x_{-1}, \ldots) \in \XX_1$, set
$\iota(\hat{x}):=(x^\ast_0,x^\ast_{-1}, \ldots) \in \XX$ so that
$x_{-k}=x^\ast_{-mk}$ for all $k \ge 0$. Set
$\XX_1^\ast:=\iota(\XX_1)$.

Since $U=g^{-1}(V) \subset V$, we have a natural lift
$\ghat^{-1}:\XX \to \XX$ of $g^{-1}$ given by $\ghat^{-1}(z_0,
z_{-1}, \ldots) := (z_{-1}, z_{-2}, \ldots)$. Note that
$\ghat^{-n}=(\ghat^{-1})^n$ embeds $\XX$ homeomorphically into
itself. Thus we can define a lift $\ghat^n:\ghat^{-n} (\XX) \to \XX$
of $g^n: g^{-n}(V) \to V$ for $n \ge 0$.

Now we claim:

\begin{proposition}\label{prop.main_step}
There exist subsets $\YY \subset \XX$ and $\YY_1 \subset \XX_1$ with
the following properties:
\begin{enumerate}[{\rm (a)}]
\item $\YY \homeo \XX$ and $\YY_1 \homeo \XX_1$ (i.e., homeomorphic).
\item Set $\YY_1^\ast := \iota(\YY_1) \subset \XX_1^\ast$. Then for all $0 \le i< m$,
$\ghat^{i}(\YY_1^\ast)$ are defined and disjoint.
\item $\YY ~\sminus~ \bigsqcup_{i=0}^{m-1} \ghat^{i}(\overline{\YY_1^\ast}) ~\homeo~ \BB_s$.
\end{enumerate}
\end{proposition}
By Proposition \ref{prop_nat_ext}, the set $\YY_1(\homeo \XX_1)$ is
homeomorphic to $\NN_{c'}$ with a compact solenoidal cone removed
where $f_{c'}$ is the straightening of $g_1$.

\parag{Remark.}
Actually we can always take $\YY=\XX$, but we can not take
$\YY_1^\ast=\XX_1^\ast$ in general. Because $(\beta(g_1),\beta(g_1),
\ldots ) \in \XX_1^\ast$ may be a fixed point of $\ghat^{-1}$ (in
the case of ``$\beta$-type" renormalizations), so we need to modify
$\XX_1^\ast$ to get the second property of the proposition.

\parag{Proof of (a) and (b).}

First we set $\YY:=\XX=\invlim(g, V)$. Then $\YY \subset \XX$ and $\YY \homeo \XX$ are trivial.

\begin{figure}[htbp]

\vspace{0cm}
\begin{center}
\includegraphics[width=.4\textwidth]{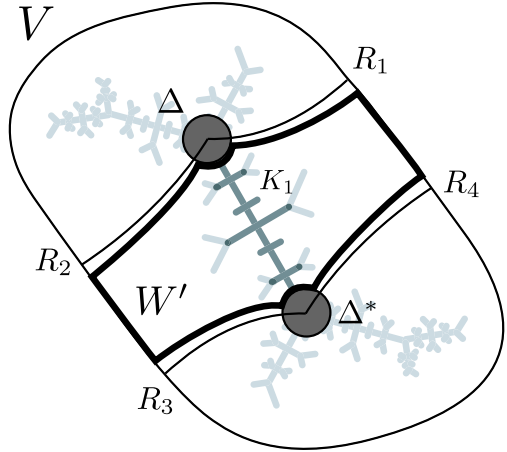}
\end{center}
\caption{The heavy curves show the boundary of $W'$.}\label{fig_W}
\end{figure}

Next we construct $\YY_1$: (In the following construction of the
topological disk $W'$, we use an idea similar to \cite[Lemmas 1.5 and
1.6]{Mil.loc.con}.) Set $\beta_1:=\beta(g_1)$ (the $\beta$-fixed
point of $g_1$) and $K_1:=K(g_1)$. Let us consider the pulled-back
external rays landing at $\beta_1$ by the straightening map $h=h_g$.
Then there are two of such rays $R_1$ and $R_2$ such that $R_1 \cup
R_2$ separates any other rays landing at $\beta_1$ and
$K_1 \sminus \skakko{\beta_1}$. Analogously, for the preimage
$\beta_1^\ast:=g_1^{-1}(\{\beta_1\}) \sminus \{\beta_1\}$, there are
two rays $R_3$ and $R_4$ landing at $\beta_1^\ast$ with the same
property. We call the rays $\{R_1,R_2,R_3,R_4\}$ the
\textit{supporting rays of} $K_1$. Let $\theta_i \in \R/\Z$ be the
angles of these $R_i$ with representatives $\theta_1<\theta_2
<\theta_3<\theta_4<\theta_1+1$. (See Figure \ref{fig_W}.)

Next we choose a sufficiently small round disk $\Delta' \subset V_1$
about $\beta_1$ so that $g_1^{-1}(\Delta')$ consists of two
topological disks $\Delta$ and $\Delta^\ast$ with $\beta_1 \in
\Delta \Subset \Delta'$ and $\beta_1^\ast \in \Delta^\ast$. We also
choose a sufficiently small $\eta>0$ such that if $t$ satisfies
$|\theta_i-t|<\eta$ for some $1 \le i \le 4$, then $R_g(t)$
intersects with either $\Delta$ or $\Delta^\ast$.

Now $V=U_g(r)$ minus the union
$$
\kakko{
R_g(\theta_1-\eta) \cup R_g(\theta_2+\eta)
\cup \overline{\Delta}
}
\sqcup
\kakko{
R_g(\theta_3-\eta) \cup R_g(\theta_4+\eta)
\cup \overline{\Delta^\ast}
}
$$
consists of three topological disks. We define $W'$ by the one
containing $K_1 \sminus \overline{\Delta}\cup \overline{\Delta^\ast}$.

Let $W_1$ denote the topological disk that is the connected
component of $W' \cap V_1$ containing the critical point of $g_1$.
Since $W_1 \subset V_1 \Subset g^{-m}(V)$, the sets $W_1$, $g(W_1)$,
$\ldots$, $g^{m-1}(W_1)$ are all defined and disjoint.

Now the inverse limit of the family $\skakko{g_1: g_1^{-n-1} (W_1)
\to g_1^{-n} (W_1)}_{n \ge 0}$, denoted by $\invlim(g_1, W_1)$, is a
proper subset of $\XX_1 = \invlim(g_1, V_1)$.

Set $\YY_1:= \invlim(g_1, W_1) \subset \XX_1$. Let us check that
$\YY_1 \homeo \XX_1$. By definition $V_1 \sminus W_1$ consists of disjoint
topological disks and does not intersects $P(g_1)$ since we take a
sufficiently small $\Delta'$. (Recall that $g$ is infinitely
renormalizable, so the $\beta$-fixed point is at a certain distance
away from the postcritical set $P(g)$. See \cite[Theorem 8.1]{Mc1}
for example. This is the only part we use the \textit{infinite}
renormalizability.) Thus $g_1: g_1^{-n-1} (V_1) \to g_1^{-n} (V_1)$
is isotopic to $g_1: g_1^{-n-1} (W_1) \to g_1^{-n} (W_1)$ for each
$n \ge 0$ and this isotopy gives a homeomorphism between the inverse
limits.

Let $\YY_1^\ast$ be the embedding of $\YY_1$ into $\XX$ by the map
$\iota$. For all $0 \le i <m$, the sets $\ghat^{i}(\YY_1^\ast)$ are
defined and disjoint since their projections $g^i(W_1)$ are defined
and disjoint. Hence we have (a) and (b) of the statement.

\parag{Proof of (c).}
Set $\BB:=\YY \sminus  \bigsqcup_{i=0}^{m-1}
\ghat^{i}(\overline{\YY_1^\ast})$. Now it is enough to show that
$\BB$ is homeomorphic to the block $\BB_s$ associated with $f_s$,
that is,
$$
\BB_s \ee \NN_s \sminus \overline{\hat{A}_s(r)} \cup
\bigsqcup_{i=0}^{m-1} \fhat_s^{i}(\overline{\hat{D}_s(1/r)}) \ee
\pi_s^{-1}(U_s(r)) \sminus \bigsqcup_{i=0}^{m-1}
\fhat_s^{i}(\overline{\hat{D}_s(1/r)}) .
$$
Here we take the same $r$ as in the construction of $V=U_g(r)$. For
later use we also set $V_s:=U_s(r)$.

\begin{figure}[htbp]
\vspace{0cm}
\begin{center}
\includegraphics[width=.7\textwidth]{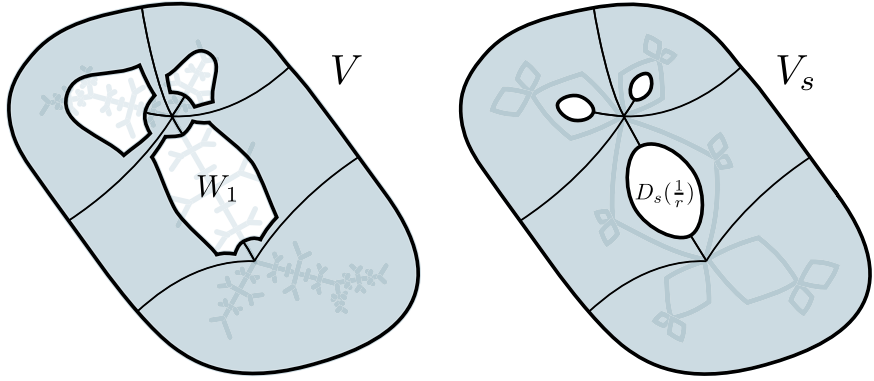}
\end{center}
\caption{The shaded region show $B$ and $B_s$ with their markings drawn in.}\label{fig_markings}
\end{figure}

We first work with the dynamics downstairs. Set $B :=V \sminus
\bigsqcup_{i=0}^{m-1} g^i(\overline{W_1})$ and mark $B$ with some
arcs given as follows (See Figure \ref{fig_markings}, left): First
join $g(\beta_1)$ and $\partial g(W_1)$ by an arc $\delta$ within
$g(\Delta)$. Since $g: W_1 \to g(W_1)$ is a branched covering, the
pull-back $g^{-1}(\delta)$ has two components in $\Delta$ and
$\Delta^\ast$. Now the markings are given by $g^{-1}(\delta)$,
$\delta$, $g(\delta), \ldots, g^{m-2}(\delta)$ and all of the
forward images of the supporting rays $\bigcup_{j=1}^4 R_j$. The
markings decompose $B$ into finitely many pieces that are all
topological disks. Note that the boundary of each piece intersects
the equipotential $E_g(r)$ and at least two external rays.

Correspondingly, set $B_s:= V_s \sminus \bigsqcup_{i=1}^{m-1}
f_s^i(\overline{D_s(1/r)})$, and complete the marking of $B_s$ by
taking all the forward images of supporting rays of $D_s$ and small
arcs from each point in the cycle $O_s$ to the corresponding
equipotentials $f_s^i(I_s(1/r))$ (Figure \ref{fig_markings}, right).
The markings also decompose $B_s$ into some pieces as in the case of
$B$.

Clearly, there is a homeomorphism $\phi$ from $B$ to $B_s$
respecting the configuration of the markings which in particular
sends the supporting external rays into the supporting external rays
without changing angles.

\parag{Lifting the homeomorphism.}
Now we claim: \textit{The map $\phi$ lifts to a homeomorphism
$\h{\phi}$ from $\pi^{-1}(B)$ in $\YY$ to $\pi_s^{-1}(B_s)$ in the
regular leaf space $\RR_s$.}

The proof requires the notion of external rays upstairs. Any
backward sequence of external rays $R_s(\theta_0) \leftarrow
R_s(\theta_{-1}) \leftarrow R_s(\theta_{-2}) \leftarrow  \cdots$
with $2 \theta_{-n}=\theta_{-n+1}$ corresponds to an arc in $\RR_s$.
Each of such arcs is parameterized by ``angles" of the form
$\h{\theta}=(\theta_0, \theta_{-1}, \ldots)$ and we denote it by
$R_s(\h{\theta})$. We define external rays $R_g(\h{\theta})$ of
$\XX$ in the same way. (Note that such angles $\{\h{\theta}\}$ and
$\sol=\invlim(f_0, \mathbb{S}^1)$ has a one-to-one correspondence.)

Recall that by construction of $W_1$, the postcritical set $P(g)$ is
contained in $\bigcup_{i=0}^{m-1} g^i(W_1)$ so $g: g^{-n-1}(B) \to
g^{-n}(B)$ is a covering map for each $n \ge 0$. Let $\Pi$ be one of
the open pieces of $B$ decomposed by the markings. Since $\Pi$ is
disjoint from the postcritical set, on each path-connected component
$\hat{\Pi}$ of $\pi^{-1}(\Pi)$ (we call it a ``plaque") the
projection $\pi|~\hat{\Pi} \to \Pi$ is a homeomorphism. Moreover,
since $\partial \Pi$ intersects with two external rays, the plaque
$\hat{\Pi}$ intersects with two external rays upstairs. Thus the
angles of these rays upstairs determine the plaques of
$\pi^{-1}(\Pi)$.

We have exactly the same situation for $B_s$. For
$\Pi_s:=\phi(\Pi)$, which is one of the compact pieces of $B_s$
disjoint from $P(f_s)$, we have a natural homeomorphic lift
$\hat{\phi} : \pi^{-1}(\Pi) \to \pi_s^{-1}(\Pi_s)$ which sends
external rays upstairs on the boundary to those without changing the
angles.

Now we have the desired homeomorphism $\hat{\phi}: \pi^{-1}(B) \to
\pi_s^{-1}(B_s)$ by gluing all of such $\hat{\phi} : \pi^{-1}(\Pi)
\to \pi_s^{-1}(\Pi_s)$ according to the angles of boundary external
rays upstairs.

\parag{Extending the homeomorphism.}
We want to extend the homeomorphism $\hat{\phi}: \pi^{-1}(B) \to
\pi_s^{-1}(B_s)$ to $\hat{\phi}: \BB \to \BB_s$. To extend
$\hat{\phi}$ to $\BB-\pi^{-1}(B)$, we first describe what this
remaining set is.

It may be easier to start with the dynamics of $f_s$. There are two
types of backward orbits which start at the cyclic Fatou components
$D_s, f_s(D_s), \ldots , f_s^{m-1}(D_s)$: One passes through $D_s$
infinitely many times, and one does only finitely many times.
Correspondingly, there are two kinds of path-connected components of
$\pi_s^{-1}(V_s \sminus B_s)$: One which is contained in the compact
solenoidal cones $\bigsqcup_{i=0}^{m-1}
\fhat_s^i(\overline{\hat{D}_s(1/r)})$, and one which is not. In
particular, the latter is a closed disk in $\BB_s$. Thus $\BB_s
\sminus \pi_s^{-1}(B_s)$ consists of such disks.

We have the same situation in the dynamics of $g$. Any
path-connected component of $\pi^{-1}(V \sminus B)$ is either
contained in the closure of $\YY_1^\ast=\iota(\invlim(g_1, W_1))$;
or not contained. The latter consists of orbits that escape from the
nest of the renormalization so it is a closed disk in $\BB$. Thus
$\BB \sminus \pi^{-1}(B)$ also consists of closed disks.

Hence it is enough to extend $\hat{\phi}$ to such ``escaping" orbits
in $\BB$. Let us choose a homeomorphic extension of $\phi$ which
maps $V$ to $V_s$ and $g^i(W_1)$ to $f_s^i(D_s(1/r))$ for all $0 \le
i<m-1$. According to the markings on $\BB$ and $\BB_s$, the
path-connected components of $\BB \sminus \pi^{-1}(B)$ and $\BB_s
\sminus \pi_s^{-1}(B)$ are labeled by the angles of external rays.
Thus there is a natural homeomorphic lift of the extended $\phi$
over those components. This gives a desired homeomorphism.

\QED(Proposition \ref{prop.main_step})

\parag{Proof of Theorem \ref{thm.blocks} (The Structure Theorem).}
One can inductively apply the argument of Proposition
\ref{prop.main_step} to each level of the renormalization
$\skakko{g_n=f_c^{p_n}|~U_n \to V_n}_{n\ge 0}$, by setting $g:=g_n$
and $g_1:=g_{n+1}$.

We first apply the proposition with $n=0$. Then we construct $W_1$
and $\BB=\BB_0$ homeomorphic to $\BB_{s_0}$. Next we apply the
proposition with $n=1$. When we take modified $g_1: U_1 \to V_1$ and
$g_2: U_2 \to V_2$ (i.e., when we replace $V_1$ by
$V_1:=U_{g_1}(r_1)$, etc.), we take a smaller $1<r_2<r_{g_2}$ so
that $U_2:=U_{g_2}(\sqrt{r_2})$ and $V_2:=U_{g_2}(r_2)$ satisfy the
original condition
$$
U_2 ~\Subset~ V_2 ~\Subset~g_1^{-m_1}(V_1)~\Subset~ U_1 ~\Subset~ V_1
$$
and the extra condition
$$
 V_2 \sqcup g_1(V_2) \sqcup \cdots \sqcup g_1^{m_1-1}(V_2) ~\Subset~W_1.
$$
As we construct $W_1 \subset V_1$, we construct $W_2 \subset V_2$ so
that $\YY_2:=\invlim(g_2, W_2)$ is homeomorphic to
$\XX_2=\invlim(g_2,V_1)$ and that $\YY_2^\ast, ~\ghat_1(\YY_2^\ast),
\ldots, \ghat_1^{m_1-1}(\YY_2^\ast)$ are defined and disjoint, where
$\YY_2^\ast$ denote the natural embedding of $\YY_2$ into
$\XX_1=\invlim(g_1,V_1)$. By the extra condition above, we have
$$
\YY_2^\ast \sqcup \ghat_1(\YY_2^\ast) \sqcup \cdots \sqcup
\ghat_1^{m_1-1}(\YY_2^\ast) ~\Subset~\YY_1.
$$
Moreover, we have a block $\BB_1' \dee \YY_1  \sminus
\bigsqcup_{i=0}^{m_1-1} \ghat_1^{i}(\overline{\YY_2^\ast})$
homeomorphic to $\BB_{s_1}$. Finally we define $\BB_1$ by the
natural embedding of $\BB_1'$ into $\XX_0 = \invlim(g_0, V_0)
\subset \NN_c$.

\begin{figure}[htbp]
\vspace{0cm}
\begin{center}
\includegraphics[width=.6\textwidth]{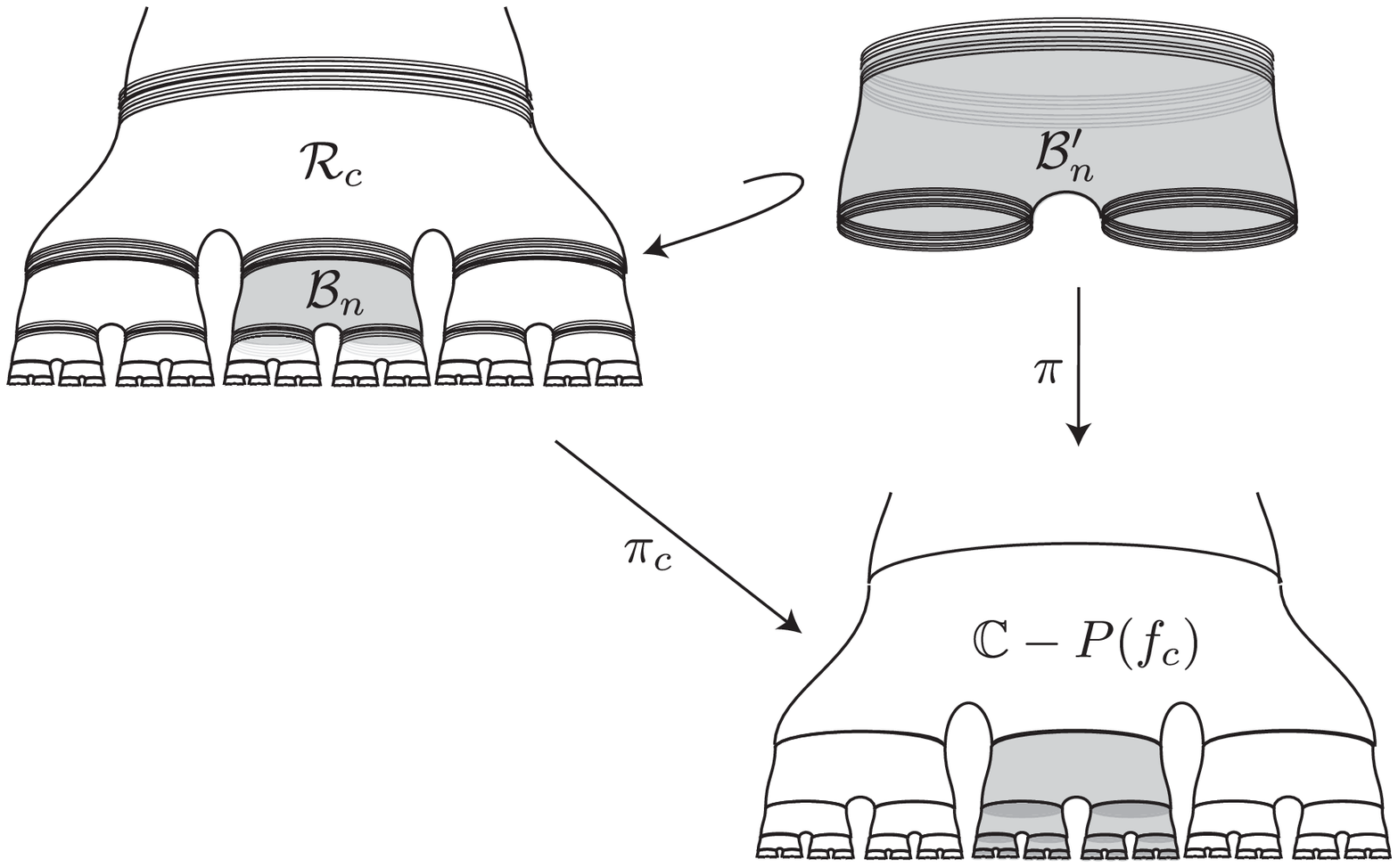}
\end{center}
\caption{A caricature of the tree structure of $\RR_c$. It comes
from a natural tree structure in the set
$\C \sminus P(f_c)$.}\label{fig_tree}
\end{figure}

Clearly the same argument works for the other levels $n \ge 2$.
Note that $\BB_n'$ constructed as above is contained in $\invlim(g_n, W_n)$.
(See Figure \ref{fig_tree}.)
So we need to iterate the natural embeddings
$$
\BB_n'
~\hookrightarrow~ \invlim(g_{n-1}, W_{n-1})
~\hookrightarrow~ \cdots
\hookrightarrow~ \invlim(g_1, W_1)
~\hookrightarrow~ \invlim(g_0, V_0)
$$
to obtain $\BB_n \subset \NN_c$.

In addition, we replace $\BB_0$ by the set $ \BB_0 \cup
\overline{\hat{A}(r_0)} \sminus \skakko{\h{\infty}} $ (where
$\hat{A}(r_0)$ is a solenoidal cone, with $r_0$ satisfying $V_0 =
U_{g_0}(r_0)$) so that $\BB_0$ covers the neighborhood of
$\hat{\infty}$. Then we have property 1 of the statement. Properties
2, 3, and 4 of the statement are clear by the construction of
blocks.

Now every backward orbit that leaves $P(f_c) \cup \skakko{\infty}$ is contained in one of
such blocks $\{\fhat_c^i(\BB_n)\}_{n,i}$. Since $f_c$ is
persistently recurrent, the set $\widehat{P(f_c)} \cup
\skakko{\h{\infty}}$ consists of all irregular points so the union
of the blocks $\{\fhat_c^i(\BB_n)\}_{n,i}$ coincide with $\RR_c$.
Thus we have 4 and 5 of the statement.

\QED

\subsection{Buildings at finite level}

To end this section we show a proposition that is important for the
arguments in the next section.

For an infinite sequence of combinatorics $\skakko{s_0, s_1,
\ldots}$, its subsequence $\skakko{s_0, s_1, \ldots, s_n}$
determines a superattracting parameter $\sigma_n$. More precisely,
for $\beta$-fixed point $\beta(g_{n+1})$ of $g_{n+1}=f_c^{p_{n+1}}$,
its forward orbit $O_{n+1}$ by $f_c$ forms a repelling periodic
point. Then its ray portrait $\mathrm{rp}(O_{n+1})$ determines a
superattracting quadratic map $f_{\sigma_n}$. It is known that it
depends only on the sub-combinatorics $\skakko{s_0, s_1, \ldots,
s_n}$ of the renormalizations.

For persistently recurrent infinitely renormalizable $f_c$ as above, we define
$$
\QQ_n \dee \RR_c \sminus
\bigsqcup_{k =n+1}^{\infty} \bigsqcup_{i=0}^{p_k-1} \fhat_c^{i}(\BB_{k}).
$$
Then we have:

\begin{proposition}\label{prop.tuning}
For $f_c$ as above, let $\QQ_n$ be the set defined as above. Then we
have a homeomorphism $h_n$ between $\QQ_n$ and $\QQ_{\sigma_n}$.
\end{proposition}

\parag{Proof.}
The proof is almost straightforward by Proposition
\ref{prop.main_step}. In fact, we can apply the same argument by
setting $g:=g_0$ and $g_1:=g_{n+1}$. \QED

\if0
An easy corollary is:
\begin{corollary}\label{cor.per_leaves}
For any periodic leaf $L$ of $\RR_{\sigma_n}$ containing a periodic
point, its preimage $h_n^{-1}(L) \cap \QQ_n$ is contained in a
periodic leaf $L'$ of $f_c$ with the same period.
\end{corollary}

\parag{Proof.}
Note that such a periodic point must be a lift of a repelling point
since infinitely renormalizable maps do not have indifferent cycles.
Thus a periodic leaf containing a periodic point is equivalent to a
leaf containing periodic external rays. By the proof of Proposition
\ref{prop.main_step}, $h_n$ sends periodic external rays to periodic
external rays. \QED
\fi

\section{Rigidity}

In this section we prove the Main Theorem of the paper which is the
following:

\begin{theorem}[Main Theorem] \label{thm.main}
Let $c$ be a non-real complex number, such that
the map $f_c$ is infinitely renormalizable with a-priori bounds. If
$h:\mc{R}_c\rightarrow \mc{R}_{c'}$ is oriented homeomorphism, then
$c$ and $c'$ belong to the same combinatorial class.
\end{theorem}

From the point of view of the parameter plane, it is known that $c$
is combinatorially rigid if and only if the Mandelbrot set is
locally connected at $c$. In view of that, our main theorem has the
following corollary.

\begin{corollary} Assume that $c$ is as in the Main Theorem and that
the Mandelbrot set is locally connected (MLC) at $c$, then $c=c'$.
\end{corollary}

In \cite{LyuI}, Lyubich proved MLC for $f_c$ with a-priori bounds
with some extra condition on combinatorics, called secondary limb
condition. In this direction, there is recent work by Jeremy Kahn
\cite{Kah1} and Kahn and Lyubich \cite{KahLyII} where they prove
a-priori bounds and MLC for infinite renormalizable parameters with
special combinatorics.

\subsection{Combinatorics of quadratic polynomials}

There are several models describing the combinatorics of quadratic
polynomials, a comprehensive text can be found in \cite{BS}, in this
paper we are going to adopt the description given by rational
laminations. Any quadratic polynomial $f_c$ with $c$ in the
Mandelbrot set, determines a relation, called the \textit{rational
lamination} of $f_c$, in $\mathbb{Q}/\mathbb{Z}$. Given $\theta$ and
$\theta'$ in $\mathbb{Q}/\mathbb{Z}$, we say that $\theta\sim
\theta'$ if the external rays $R_{\theta}$ and $R_{\theta'}$ land at
the same point in the Julia set $J(f_c)$. Jan Kiwi gave a set of
properties which guarantee that if a given relation in
$\mathbb{Q}/\mathbb{Z}$ satisfies these properties then the relation
is a rational lamination of some polynomial $P$, the interested
reader can consult \cite{kiwi}. For us, the most relevant property
of rational laminations is the following:

\begin{lemma}\label{comb.rat.lam} Let $R$ and $R'$ be two rational
laminations, assume that there is $\theta\in \mathbb{Q}/\mathbb{Z}$
such that each class in $R'$ is obtained by rotating a class in $R$
by angle $\theta$. Then $\theta=0\, \mod(1)$.\end{lemma}

Let us call a leaf $L$ in $\mc{R}_c$ \textit{repelling} if it
contains a repelling periodic point of $\h{f}_c$. Clearly, every
repelling leaf is invariant under some iterate of $\h{f}_c$, the
converse is not true in general, because in the presence of
parabolic point there are invariant leaves without periodic points.
In the case of the dyadic solenoid $\sol$ if a leaf $L\subset \sol$
is invariant under some iterate of $\h{f}_0$, then $L$ is repelling.
The fact that all periodic leaves in $\sol$ are repelling allow us
to lift combinatorial properties of periodic points in $J(f_c)$ to
repelling leaves in $\mc{R}_c$.

More precisely, let $L$ be a repelling leaf in $\mc{R}_c$ and let
$\mc{S}_r$ some solenoidal equipotential, the intersection $L\cap
\mc{S}_r$ consists of some leaves in $\mc{S}_r$ under the canonical
identification, it turns out, every such leaf is repelling in $\sol$
under $\h{f}_0$. Moreover, the pullback to $L$ of each of these
periodic points is precisely the intersection of a periodic
solenoidal external ray landing at the periodic point of $L$.

In the dynamical plane, if $p$ is a periodic point in the Julia set
$J(f_c)$ then $p$ is the landing point of external rays which are
periodic under $f_c$, see \cite{Mext}, if the periodic lift $\h{p}$
belongs to the regular part, then there are periodic solenoidal
external rays landing at $\h{p}$ in $L(\h{p})$, each of these
solenoidal external rays will intersect a leaf of a solenoidal
equipotential. As a consequence we have:

\begin{lemma}\label{lem.com.leaf} Two leaves
$\mc{S}_r$, coming from periodic leaves in $\sol$, belong to the
same leaf $L$ in $\mc{R}_c$ if and only if they intersect periodic
solenoidal external rays landing in the same point in
$\fiber(J(f_c))\cup \mc{R}_c$.
\end{lemma}

We will see that, for quadratic polynomials with a-priori bounds,
repelling leaves have topological relevance. Such was the approach
in \cite{cab2} (see also \cite{cab}) to prove rigidity for
hyperbolic maps and complex semi-hyperbolic. We can resume the main
results in \cite{cab2} with the following theorem:

\begin{theorem}\label{comb.class} Let $h:\mc{N}_c\rightarrow \mc{N}_{c'}$ be a
homeomorphism between natural extensions, such that:
\begin{enumerate}
\item $h(\h{\infty})=\h{\infty}$; and \item $h$ sends repelling leaves
into repelling leaves.
\end{enumerate}Then $f_c$ and $f_{c'}$ belong to the same combinatorial class.
\end{theorem}

The proof of this theorem is decomposed in three statements;
Lemma~\ref{first.lemma} whose proof can be found in \cite{cab2},
Proposition~\ref{Kwapisz} due to Jaroslaw Kwapisz \cite{Kwap}, and
Lemma~\ref{last.lemma}. The first starts by noting that the
foliation of the solenoidal cone by solenoidal equipotentials
defines a local base of neighborhoods at $\h{\infty}$ in $\mc{N}_c$.
Hence, given a homeomorphism $h$ as in Theorem~\ref{comb.class}, we
can find a solenoidal equipotential $\mc{S}_r$ whose image lies
between two solenoidal equipotentials. Recall that a solenoidal
equipotential $\mc{S}_r$ has associated a canonical homeomorphism
$\phi_R:\mc{S}_r\rightarrow \sol$, moreover,
$\phi_{R^2}\circ\h{f}_c\circ\phi_R^{-1}=\h{f}_0$. Hence, we are in
the following situation:

\begin{lemma}\label{first.lemma}Let $e:\sol \rightarrow \sol \times (0,1)$ be
a topological embedding, then there is a map $e'$ isotopic to $e$
such that $e'(\sol)=\sol\times {1/2}$.
\end{lemma}

We can pull back the isotopy in this lemma, to an isotopy defined on
$\mc{S}_r$ which extends to an isotopy defined on a neighborhood of
$\mc{S}_r$. Hence, we can find a homeomorphism $h'$, isotopic to
$h$, that sends homeomorphically a solenoidal equipotential into a
solenoidal equipotential. With the canonical homeomorphism of
solenoidal equipotentials to $\sol$, $h'$ induces a self
homeomorphism of the dyadic solenoid $\sol$. Now, as described by
Kwapisz in \cite{Kwap}, each homotopic class of homeomorphisms of
$\sol$ is uniquely represented by a map with a special form:

\begin{proposition}[Kwapisz]\label{Kwapisz} Let $\phi:\sol \rightarrow \sol$ be a
homeomorphism of the dyadic solenoid, then there exist $n$ and an
element $\tau\in \sol$ such that $\phi$ is isotopic to $\h{z}\mapsto
\tau \h{f}^n_0(\h{z})$.
\end{proposition}

The number $n$ is uniquely determined by the homotopic class of
$h'$, so if we post-compose $h'$ with $f^{-n}_{c'}$,
Proposition~\ref{Kwapisz} implies that we can find a new
homeomorphism from $\mc{N}_c$ to $\mc{N}_{c'}$ sending one
solenoidal equipotential into a solenoidal equipotential, such that
under the canonical identification, the map between these solenoidal
equipotentials is just the left translation by $\tau$ of the dyadic
solenoid $\sol$.

All isotopies above, and the map $\h{f}_{c'}$, send repelling leaves
into repelling leaves, so our new homeomorphism will also send
repelling leaves into repelling leaves. By the previous lemmas, if
$h$ is a homeomorphism like in Theorem~\ref{comb.class}, then we can
assume that $h$ sends a solenoidal equipotential $\mc{S}_r$
homeomorphically into a solenoidal equipotential and, that under
canonical isomorphisms, the map $h$ restricted to $\mc{S}_r$ is just
a translation $\tau$ by an element in $\sol$. Now the combinatorial
information of $f_c$ give us more restrictions on the isotopy class
of $h$:

\begin{lemma}\label{last.lemma} Assume $h$ is a homeomorphism like
in Theorem~\ref{comb.class}, then the induced translation $\tau$ in
Proposition~\ref{Kwapisz} is homotopic to the identity.
\end{lemma}

\parag{Proof.} Let us consider the restriction of $h$ to the solenoidal
equipotential $\mc{S}_r$ such that $h(\mc{S}_r)$ is also a
solenoidal equipotential, under canonical homeomorphisms the map
$H=h|\mc{S}_r$ is a map from $\sol$ into itself. We assume that $H$
has the form $\h{z}\mapsto \tau \h{z}$. By Lemma~\ref{lem.com.leaf},
$h|\mc{S}_r$ sends repelling leaves into repelling leaves.

Let $L$ be a periodic leaf in $\sol$ with $\h{\theta}$ the periodic
point in $L$, let $\h{\theta}'$ be the periodic point in $H(L)$. By
sliding $\sol$ along $h|\mc{S}_r(L)$ to send $H(\h{\theta})$ to
$\h{\theta}'$, this operation induces a new map $H'$ in the isotopy
class of $H$, which satisfies
$H(\h{\theta})=\tau'\h{\theta}=\h{\theta}$, since $\h{\theta}$ and
$\h{\theta}'$ are periodic in $\sol$, $\tau'$ must be periodic as
well. Hence, the map $H'$ leaves the set of periodic points in
$\sol$ invariant.

Now, periodic points in $\sol$ are determined by the first
coordinate. The translation $\tau$ induces a rotation in the set of
periodic angles which extends to a rotation on the rational
lamination. By Lemma~\ref{comb.rat.lam} this implies that the
rational laminations are the same, and that the translation $\tau'$
is the identity, by construction $\tau'$ is isotopic to $\tau$. \QED

\parag{Proof.}[Proof of Theorem \ref{comb.class}] As a consequence of
the previous Lemma the rational laminations of $f_c$ and $f_{c'}$
are the same. This implies that $c$ and $c'$ belong to the same
combinatorial class. \QED

\subsection{Ends of the regular part}

A path $\gamma:[0,\infty)\rightarrow \mc{R}_c$ is said to
\textit{escape to infinity} if it leaves every compact set $K\subset
\mc{R}_c$. we define an \textit{end} of $\mc{R}_c$ to be an
equivalence class of paths escaping to infinity. Let $\gamma$ and
$\sigma$ two paths escaping to infinity, we say that $\gamma$ and
$\sigma$ access the same end if for every compact set $K\subset
\mc{R}_c$, the paths $\gamma$ and $\sigma$ eventually belong to the
same connected component of $\mc{R}_c\setminus K$. Consider the set
$End(\mc{R}_c)$ consisting of $\mc{R}_c$ union with the abstract set
of ends.

Let $f_c$ be an infinitely renormalizable quadratic polynomial with
a-priori bounds, by Theorem \ref{KL.loc.comp} the regular part
$\mc{R}_c$ is locally compact and then $End(\mc{R}_c)$ is a compact
set, which we will call the \textit{end compactification} of
$\mc{R}_c$.

\begin{proposition}\label{irr.end} Let $f_c$ be an
infinite renormalizable quadratic polynomial with a-priori bounds,
then $End(\mc{R}_c)$ is homeomorphic to $\mc{N}_c$.
\end{proposition}

\parag{Proof.} We will show that there exist a bijection $\Phi$ between the
set of irregular points and the set of ends. Let $\h{i}$ be an
irregular point in $\mc{N}_c$, let $i_0=\pi(\h{i})$ and take any
$z_0\in \C\setminus \omega(c)$. Since $\omega(c)$ is a Cantor set,
there is a path $\sigma$ be a path connecting $z_0$ with $i_0$ which
intersects $\omega(c)$ only at $i_0$. We can lift the path $\sigma$
to $\mc{N}_c$ to a path $\h{\sigma}$ from a point in the fiber of
$z_0$ connecting to $\h{i}$. By construction, the path $\h{\sigma}$
intersects $\mc{I}_c$ at $\h{i}$, then the restriction of
$\h{\sigma}$ to $\mc{R}_c$ is a path escaping to infinity. Let
$\Phi(\h{i})=[\h{\sigma}]$, where $[\h{\sigma}]$ is the end
represented by $\h{\sigma}$. Now we check that $\Phi$ is well
defined, let $\h{\sigma}$ and $\h{\sigma'}$ be two paths in
$\mc{N}_c$ intersecting the irregular set only at the end point
$\h{i}$. These paths do not need to start at the same point or
belong to the same leaf. Let $L$ be the leaf containing
$\sigma([0,1))$ in $\mc{R}_c$. Since every leaf is dense in
$\mc{R}_c$ and is simply connected, we can construct a family of
paths $\h{\sigma_n}$ in $L$, ending at $\h{i}$ and such that
$\h{\sigma}_n\rightarrow \h{\sigma}'$ pointwise. Let $K$ be any
compact set in $\mc{R}_c$, and $U$ be a connected component of
$\mc{R}_c\setminus K$ which eventually contains $\h{\sigma}'$. Since
$U$ is open, there is a $N$ such that $\h{\sigma}_N$ also eventually
belongs to $U$, but $\h{\sigma}$ and $\h{\sigma}_N$ belong to the
same path connected component (same leaf), thus $\h{\sigma}$ must
also be eventually contained in $U$.

To see that $\Phi$ is injective, let $\h{i}$ and $\h{i}'$ be two
irregular points, since the projection $\pi$ is a homeomorphism
between the set of irregular points and $\omega(c)$ we have
$\pi(\h{i})\neq \pi(\h{i}')$, and any two paths $\phi$ and $\phi'$
escaping to $\h{i}$ and $\h{i}'$ respectively, must eventually
belong to different components of some level of renormalization.

Finally, let us prove that $\Phi$ is surjective. Let $e$ be an end
of $\mc{R}_c$, and consider $\phi$ a path escaping to $e$. Let $D_r$
be a closed ball containing $J(f)$. For each level of
renormalization $n$, let $Q_n$ be a family of disjoint open
neighborhoods of the little Julia sets of level $n$, if these Julia
set touch, we can shrink the domains a little to make them disjoints
as in the proof of Theorem \ref{thm.blocks}. Let $W_n$ be the union
of the domains in $Q_n$. Then $K_n=B_r\setminus Q_n$ is a compact
set in $\C\setminus \omega_c$. Thus $\fiber(K_n)$ is compact in
$\mc{R}_c$, by definition the path $\phi$ must eventually escape
$\fiber(K_n)$. It follows that the projection $\pi(\phi)$ eventually
belongs either to a neighborhood of infinity, and then $\phi$
escapes to $\h{\infty}$, or to a domain in $Q_n$, say $V_n$, by the
disjoint property of the sets in $Q_n$, it is clear that $V_{n+1}$
is contained in $V_n$. By construction, the domains $\{V_n\}$ shrink
to a point $i_0$ in $\omega(c)$.  This process can be repeated for
every coordinate of $\phi$ to get a sequence of points $\{i_n\}$ in
$\omega(c)$ which are the coordinates of a point $\h{i}$ in
$\widehat{\omega(c)}$. Since $f$ is persistently recurrent $\h{i}$ is irregular. \QED\\

On the remaining part of the paper $h$ will denote a homeomorphism
of the regular parts of two infinite renormalizable quadratic
polynomials $f_c$ and $f_{c'}$ with a-priori bounds.

\begin{corollary}\label{fix.infty} The map $h$ admits an extension to a homeomorphism
$\tilde{h}:\mc{N}_c\rightarrow \mc{N}_{c'}$ of the regular
extensions. Moreover, $\tilde{h}(\h{\infty})=\h{\infty}$.
\end{corollary}
\parag{Proof.} By Proposition \ref{irr.end} the map $h$ extends to the natural
extensions sending irregular points to irregular points, and by
Lemma~\ref{feig.irreg} the point $\h{\infty}$ is the only isolated
irregular point, hence $h(\h{\infty})=\h{\infty}$. \QED

\subsection{Topology of Periodic leaves}

Since leaves are path connected components of $\mc{R}_c$, given a
leaf $L\subset \mc{R}_c$ we can consider how many access to
$\h{\infty}$ the leave has. That is, the number of path components
of $L\setminus K$ that are connected to $\h{\infty}$ in $\mc{N}_c$,
for a suitable large compact set $K\subset \mc{R}_c$. Note that a
leaf has access to points in $\omega(c)$ if and only if intersects
infinitely many levels in the tree structure of $\mc{R}_c$. However,
this is not the case for repelling leaves:

\begin{lemma} Let $L$ be a repelling leaf, then $L$ there is a level
$n$ such that $L\subset \mc{Q}_n$. In this case, $L$ has access
only to $\h{\infty}$. \end{lemma}

\parag{Proof.}Let $\h{p}$ be the periodic point in $L$ and let $p=\pi(p)$.
Since $f_c$ is infinite renormalizable, $p$ is repelling, and
therefore it must belong to the Julia set $J(f_c)$, moreover, the
inverse of the classical K\"{o}nigs linearization coordinate
around $p$ provides a global uniformization coordinate for $L$.
From this uniformization it follows that a point $\h{z}$ in
$\mc{R}_c$ belongs to $L$ only if the coordinates of $\h{z}$
converge to the cycle of $p$.

Since the intersection of the renormalization domains is just the
postcritical set, we can find a level $n+1$ of the renormalization
such that the orbit of the renormalization domains of level $n+1$ is
outside a neighborhood of the cycle of $p$. By this choice, no point
in $L$ can intersect the level $n+1$ of the tree structure of
$\mc{R}_c$. The statement of the lemma now follows. \QED

When $f_c$ is superattracting, every leaf $L$ invariant under some
iterate of $\h{f}$ must contains a repelling periodic point and
hence $L$ is repelling. In this case, there are no critical points
in the Julia set $J(f_c)$ so the fiber $\fiber(J(f_c))$ is compact.
If $p$ is a periodic point in $J(f_c)$. Let $\h{p}$ be invariant
lift of $p$ in $\mc{R}_c$, and $L(\h{p})$ the leaf containing
$\h{p}$. From \cite{cab2}, we have the following:

\begin{proposition}\label{per.sup} The number of access of $L$ to $\h{\infty}$ is
equal to the number of external rays landing at $p$. Moreover, if
$L$ is a leaf which has at least three access to infinity, then $L$
must be repelling.
\end{proposition}
Let us remark that in the superattracting case,
Proposition~\ref{irr.end} also holds, however, repelling leaves may
have access to other irregular points. Nevertheless, if some
repelling leaf $L$ has at least three access to $\h{\infty}$ then by
Proposition~\ref{per.sup}, the corresponding periodic point $p$ has
at least three external rays landing at $p$. This situation
only can happen if the imaginary part of $c$ is not 0.\\

Let us now go back to the case were $f_c$ is infinite
renormalizable with a-priori bounds:

\begin{lemma}\label{im.per} Let $f_c$ be infinite renormalizable with a-priori
bounds, and let $L\subset\mc{R}_c$ be a leaf which has at access
only to $\h{\infty}$, and such that the number of access to infinity
is at least 3, then $L$ must be a repelling leaf. Moreover, this
implies that $Im(c)\neq 0$.
\end{lemma}

\parag{Proof.} Since the only access to infinity of $L$ is
$\h{\infty}$, there is a level $n$ such that $L\subset \mc{Q}_n$. By
Corollary~\ref{fix.infty} the map $h$ extends to the natural
extensions and $\h{\infty}$, so the image $h(L)$ is also a leave
with the same number of access to $\h{\infty}$. Regarding $L$ as a
subset of $\mc{Q}_n$, the leaf $L$ has at least 3 access to $\infty$
in $\mc{Q}_n$ by Proposition~\ref{per.sup} the leaf $L$ must be
repelling in $\mc{Q}_n$ under dynamics of $\h{f}_{s_n}$, by the
block homeomorphism in the proof of Theorem~\ref{prop.tuning}, this
implies that $L$ itself must be repelling under dynamics of
$\h{f}_c$. \QED\\

Now we are ready to prove the Main Theorem:

\parag{Proof of Main Theorem.} By Corollary \ref{fix.infty}, the map
extends to a homemorphisms of natural extensions $\tilde{h}$ with
$\tilde{h}(\h{\infty})$. Since $Im(c)\neq 0$ then there exist a
repelling leaf $L$ in $\mc{N}_c$ such that $L$ has at least three
access to $\h{\infty}$. This is a topological property, so $h(L)$ is
also a leaf with at least 3 access to $\h{\infty}$. By Lemma
\ref{im.per} $h(L)$ is also repelling and moreover $Im(c')\neq 0$.
In this way, $\tilde{h}$ sends a repelling leaf into a repelling
leaf. By an isotopy argument similar to the one used in the proof of
Lemma \ref{last.lemma}, we can see that this implies that $h$ sends
repelling leaves into repelling leaves. Hence, $\tilde{h}$ satisfies
the conditions of Theorem \ref{comb.class}, which implies that $f_c$
and $f_{c'}$ belong to the same combinatorial class. \QED

\bibliographystyle{amsplain}
\bibliography{workbib}

\begin{tabular}{ll}
\hspace{7cm} &\quad Carlos Cabrera  \\
 &\quad Institute of Mathematics of the\\
 &\quad Polish Academy of Sciences\\
   &\quad ul. \'{S}niadeckich 8\\
   &\quad 00-956 Warszawa, Poland \\

\end{tabular}

\bigskip

\begin{tabular}{ll}
\hspace{7cm} &\quad Tomoki Kawahira \\
 &\quad Graduate School of Mathematics \\
 &\quad Nagoya University \\
   &\quad   Nagoya 464-8602, Japan \\

\end{tabular}

\end{document}